\documentclass[11pt,a4paper]{article}
\usepackage{amsmath,amssymb}
\usepackage{latexsym}
\usepackage{mathrsfs}
\usepackage{amsfonts}
\usepackage{hyperref}
\usepackage{amsthm}
\usepackage{appendix}
\usepackage[square,comma,sort&compress,numbers]{natbib}
\usepackage{color}
\usepackage{subcaption}
\usepackage{wrapfig}
\usepackage{tikz}
\usepackage{graphicx}

\hfuzz=\maxdimen
\theoremstyle{definition}
\newtheorem{lemma}{Lemma}[section]
\newtheorem{definition}[lemma]{Definition}
\newtheorem{theorem}[lemma]{Theorem}

\newtheorem{proposition}[lemma]{Proposition}
\newtheorem{remark}{Remark}
\newtheorem{convention}[lemma]{Convention}
\newtheorem{question}[lemma]{Question}

\DeclareFixedFont{\Acknowledgment}{OT1}{cmr}{bx}{n}{14pt}
\begin{document}

\title{\bf Prescribing discrete Gaussian curvature on polyhedral surfaces}
\author{Xu Xu, Chao Zheng}
\maketitle

\begin{abstract}
Vertex scaling of piecewise linear metrics on surfaces introduced by Luo \cite{Luo1}
is a straightforward discretization of smooth conformal structures on surfaces.
Combinatorial $\alpha$-curvature for vertex scaling of piecewise linear metrics on surfaces
 is a discretization of Gaussian curvature on surfaces.
In this paper, we investigate the prescribing combinatorial $\alpha$-curvature problem on polyhedral surfaces.
Using Gu-Luo-Sun-Wu's discrete conformal theory \cite{GLSW} for piecewise linear metrics on surfaces
and variational principles with constraints,
we prove some Kazdan-Warner type theorems for prescribing combinatorial $\alpha$-curvature problem,
which generalize the results obtained in \cite{GLSW, Xu3} on prescribing combinatorial curvatures on surfaces.
Gu-Luo-Sun-Wu \cite{GLSW} conjectured that one can prove Kazdan-Warner's theorems in \cite{KW1,KW2} via approximating
smooth surfaces by polyhedral surfaces.
This paper takes the first step in this direction.
\end{abstract}



\section{Introduction}
\subsection{Motivation}
A  classical problem in differential geometry is prescribing Gaussian curvature on closed surfaces,
which asks the follow question.
\\
\\
\textbf{Prescribing Gaussian curvature problem.} Let $S$ be a smooth connected closed surface.
What functions could be realized as Gaussian curvatures of Riemannian metrics on $S$?
\\
\\
This problem has been  extensively studied by lots of mathematicians
by deforming Riemannian metric in its conformal structure.
See, for instance, \cite{Berger, CY, Moser, KW1, KW2} and others.

Polyhedral metric on surfaces is a discrete analogue of the Riemannian metric on surfaces,
which is a constant curvature ($1$, $0$ or $-1$) metric with cone singularities.
The classical discrete Gaussian curvature for a polyhedral metric on a surface is
the corresponding discrete analogue of the smooth Gaussian curvature,
which is defined to be the angle defect at the cone points.
Prescribing discrete Gaussian curvature problem
on polyhedral surface asks the following question parallelling to the classical prescribing Gaussian curvature problem on smooth surfaces.
\\
\\
\textbf{Prescribing discerete Gaussian curvature problem.}
Suppose $S$ is a connected closed surface and  $V$ is a finite nonempty subset of $S$.
What functions defined on $V$ could be realized as discrete Gaussian curvatures of polyhedral metrics on $(S, V)$?
\\
\\
In discrete conformal geometry,
parallelling to the smooth case, the prescribing discrete Gaussian curvature problem on surfaces
is usually studied by deforming a polyhedral metric in its discrete conformal structure,
which is a discrete analogue of  the smooth conformal structure on surfaces defining polyhedral metrics
by functions defined on the cone points.
The discrete conformal structures on surfaces
that have been extensively studied in the history include Thurston's circle packings \cite{T1},  Luo's vertex scaling \cite{Luo1}  and others.
For Thurston's circle packings on surfaces, the solution of prescribing discrete Gaussian curvature problem
gives rise to the famous Koebe-Andreev-Thurston Theorem, which plays a fundamental role in the study
of geometry and topology of 3-dimensional manifolds \cite{T1}.
For vertex scaling of piecewise linear and piecewise hyperbolic metrics on surfaces,
the recent perfect solution of prescribing discrete Gaussian curvature problem by Gu-Luo-Sun-Wu \cite{GLSW}
and Gu-Guo-Luo-Sun-Wu \cite{GGLSW}
provides a new constructive proof of the classical uniformization theorem on closed surfaces with genus $g\geq 1$,
which is computable and has lots of applications \cite{DGL,GY, Luo4, ZG, GLY, SWGL}.
In the framework of vertex scaling,
Springborn \cite{springborn} recently proved that the discrete uniformization theorem on the sphere
is equivalent to Rivin's realization theorem for ideal hyperbolic polyhedra \cite{Rivin2}.

However, the classical discrete Gaussian curvature is not a proper discretization of the smooth Gaussian curvature on surfaces.
For example, it is scaling invariant, which is different from the transformation of smooth Gaussian curvatures
under scaling of Riemannian metrics,
and does not approximate the smooth Gaussian curvature pointwisely as the triangulations become finer and finer \cite{GX2}.
In \cite{Xu3}, the first author introduced combinatorial $\alpha$-curvature for vertex scaling
of piecewise linear metrics (PL metrics for short in the following)
on surfaces to avoid the disadvantages of classical discrete Gaussian curvature alluded to above,
where the rigidity and Yamabe problem for combinatorial $\alpha$-curvature were also studied.
In this paper, we study the prescribing combinatorial $\alpha$-curvature problem for PL metrics
on surfaces and prove some Kazdan-Warner type theorems.
It is conjectured by Gu-Luo-Sun-Wu \cite{GLSW} that
one can prove Kazdan-Warner's theorems in \cite{KW1,KW2} by approximating
smooth surfaces by polyhedral surfaces involving vertex scaling.
This paper takes the first step in this direction.

\subsection{Statements of main results}

Suppose $S$ is a connected closed surface and $V$ is a finite non-empty subset of $S$, $(S, V)$ is called a marked surface.
A PL metric $d$ on the marked surface $(S,V)$ is a flat cone metric on $S$ with conic singularities contained in $V$.
A triangulation $\mathcal{T}$ of the marked surface $(S, V)$ is a triangulation of $S$ with the vertex set equal to $V$.
Denote a vertex, an edge and a face in the triangulation $\mathcal{T}$ by $i, \{ij\}, \{ijk\}$ respectively and denote the sets of edges and faces of $\mathcal{T}$ by $E_\mathcal{T}$ and $F_\mathcal{T}$ respectively.
A triangulation $\mathcal{T}$ of $(S,V,d)$ is geometric if
every edge in $E_\mathcal{T}$ is a geodesic in the PL metric $d$.
A PL metric $d$ on a triangulated surface $(S,V, \mathcal{T})$ with $\mathcal{T}$ geometric
defines a map $l: E_\mathcal{T}\rightarrow (0, +\infty)$ such that $l_{ij}, l_{ik}, l_{jk}$ satisfy the triangle inequalities for any triangle $\{ijk\}\in F_\mathcal{T}$.
The map $l: E_\mathcal{T}\rightarrow (0, +\infty)$ is called as a discrete metric and the logarithm of the discrete metric $l$,
$\lambda_{ij}=2\log l_{ij}$, is called as the logarithmic length.
One can also obtain discrete PL metrics on a triangulated surface $(S,V, \mathcal{T})$ by gluing triangles in 2-dimensional Euclidean space isometrically along the edges in pair, which gives rise to PL metrics on the marked surface $(S, V)$.
Note that a PL metric $d$ on a marked surface $(S,V)$ is intrinsic in the sense
that it is independent of the geometric triangulations of $(S,V, d)$.
Every PL metric $d$ on $(S,V)$ has a Delaunay triangulation $\mathcal{T}$ of $(S,V, d)$
such that each triangle in $\mathcal{T}$ is Euclidean and the sum of two angles facing each edge is at most $\pi$.
See \cite{Aurenhammer,BS,GLSW,Rivin} for further discussion on Delaunay triangulation of polyhedral surfaces.

Suppose $(S,V,d)$ is a marked surface with a PL metric $d$ and $\mathcal{T}$ is a geometric triangulation of $(S,V,d)$.
The classical discrete Gaussian curvature $K: V\rightarrow (-\infty, 2\pi)$ for $(S,V,d)$ is defined as
\begin{equation}\label{K_i}
K_i=2\pi-\sum_{\{ijk\}\in F_\mathcal{T}}\theta_i^{jk}
\end{equation}
with summation taken over all the triangles at $i\in V$ and $\theta_i^{jk}$ being the inner angle of the triangle $\{ijk\}\in F_\mathcal{T}$ at the vertex $i$.
Sometimes we call $K$ as combinatorial curvature for simplicity.
Note that the classical discrete Gaussian curvature $K$ is intrinsic in the sense that
it is independent of the geometric triangulations of $(S,V,d)$.
The classical combinatorial curvature $K$ for PL metrics on $(S,V)$
satisfies the following discrete Gauss-Bonnet formula (\cite{Chow-Luo}, Proposition 3.1)
\begin{equation}\label{discrete Gauss-Bonnet formula}
\sum_{i\in V}K_i=2\pi\chi(S).
\end{equation}
The discrete Gauss-Bonnet formula (\ref{discrete Gauss-Bonnet formula}) provides a necessary condition for
a function to be the discrete Gaussian curvature of some PL metric on $(S, V)$.
The classical prescribing discrete Gaussian curvature problem on surfaces can be taken
as a converse problem to the discrete Gauss-Bonnet formula (\ref{discrete Gauss-Bonnet formula}), which asks the following question.
\begin{question}
Is the discrete Gauss-Bonnet formula (\ref{discrete Gauss-Bonnet formula})
sufficient condition for a function $K: V\rightarrow (-\infty, 2\pi)$ to be the discrete Gaussian curvature of some PL metric on $(S, V)$?
\end{question}

The prescribing discrete Gaussian curvature problem is usually studied in the frame of discrete conformality
of polyhedral metrics on surfaces. See, for instance, \cite{T1,GLSW,GGLSW, Luo1} and others.
Vertex scaling of PL metrics on surfaces introduced by Luo \cite{Luo1}
is a straightforward discrete analogue of the conformal transformation in Riemannian geometry.
\begin{definition}[\cite{Luo1}]\label{luo's discrete conformal}
Suppose $l,\widetilde{l}: E_{\mathcal{T}}\rightarrow (0, +\infty)$ are two discrete PL metrics on a triangulated surface $(S,V, \mathcal{T})$.
$\widetilde{l}$ is called as a vertex scaling of $l$ if
there exists a function $u: V\rightarrow \mathbb{R}$ such that
$$\widetilde{l}_{ij}=l_{ij}\exp\left(\frac{u_i+u_j}{2}\right)$$
for any edge $\{ij\}\in E_\mathcal{T}$.
\end{definition}

Among other things, Luo \cite{Luo1} further established a variational principle for vertex scaling of PL metrics on surfaces
and proved that there
exists some combinatorial obstructions for the existence of PL metric with constant discrete Gaussian curvature $K$
in a discrete conformal class on a triangulated surface $(S, V, \mathcal{T})$ in the sense of Definition \ref{luo's discrete conformal}.
This implies that there exist some combinatorial obstructions for the solvability of the prescribing discrete Gaussian curvature problem
in a discrete conformal class on a triangulated surface $(S, V, \mathcal{T})$  in the sense of Definition \ref{luo's discrete conformal}.
To overcome this difficulty, based on Bobenko-Pinkall-Springborn's important work \cite{BPS},
Gu-Luo-Sun-Wu \cite{GLSW} introduced the following new definition of discrete conformality of PL metrics on marked surfaces,
which allows the triangulation of the marked surface to be changed under the Delaunay condition.

\begin{definition}[\cite{GLSW}, Definition 1.1]\label{discrete conformal equivalent}
Two PL metrics $d,d'$ on a marked surface $(S,V)$ are discrete conformal
if there exist a sequence of PL metrics $d_1=d,d_2,...,d_m=d'$ on $(S,V)$ and triangulations $\mathcal{T}_1,...,\mathcal{T}_m$ of $(S,V)$ satisfying
\begin{description}
  \item[(a)] (Delaunay condition) each $\mathcal{T}_i$ is Delaunay in $d_i$,
  \item[(b)] (Vertex scaling condition) if $\mathcal{T}_i=\mathcal{T}_{i+1}$, there exists a function $u: V\rightarrow \mathbb{R}$, called a conformal factor, so that if $e$ is an edge in $\mathcal{T}_i$ with end points $v$ and $v'$, then the lengths $l_{d_i}(e)$ and $l_{d_{i+1}}(e)$ of $e$ in metrics $d_i$ and $d_{i+1}$ are related by $$l_{d_{i+1}}(e)=l_{d_i}(e)\exp\left(\frac{u(v)+u(v')}{2}\right),$$
  \item[(c)] if $\mathcal{T}_i\neq\mathcal{T}_{i+1}$, then $(S,d_i)$ is isometric to $(S,d_{i+1})$ by an isometry homotopic to the identity in $(S,V)$.
\end{description}
\end{definition}
Definition \ref{discrete conformal equivalent} defines an equivalence relationship for PL metrics on a marked surface $(S,V)$.
The equivalence class of a PL metric $d$ on $(S,V)$ is called as the discrete conformal class of $d$ and
denoted by $\mathcal{D}(d)$.
Using the new discrete conformality in Definition \ref{discrete conformal equivalent},
Gu-Luo-Sun-Wu \cite{GLSW} perfectly solved the prescribing discrete Gaussian curvature problem for PL metrics on closed surfaces
in the following well-known theorem, which shows that the discrete Gauss-Bonnet formula (\ref{discrete Gauss-Bonnet formula})
is a necessary and sufficient condition for a function $\overline{K}: V\rightarrow (-\infty, 2\pi)$
to be the classical discrete Gaussian curvature of some PL metric on $(S, V)$.
\begin{theorem}[\cite{GLSW} Theorem 1.2]\label{Euclidean discrete uniformization}
Suppose $(S, V)$ is a closed connected marked surface and $d$ is a PL metric on $(S, V)$.
Then for any $\overline{K}: V\rightarrow (-\infty, 2\pi)$ with $\sum_{v\in V}\overline{K}(v)=2\pi\chi(M)$,
there exists a PL metric $\overline{d}$, unique up to scaling and isometry homotopic to the identity
 on $(S, V)$, such that $\overline{d}$ is discrete conformal to $d$ and the discrete curvature
 of $\overline{d}$ is $\overline{K}$.
\end{theorem}

The classical discrete Gaussian curvature defined by (\ref{K_i})
is not a proper discretization of the smooth Gaussian curvature on surfaces,
which is supported by the discussions in \cite{BPS, GX2}.
To overcome the disadvantages of the classical discrete Gaussian curvature,
Ge and the first author \cite{GX2} introduced the combinatorial $\alpha$-curvature for Thurston's Euclidean circle packing metrics on surfaces.
After that, there are lots of research activities on combinatorial $\alpha$-curvature on surfaces and $3$-manifolds.
See, for instance, \cite{DG,GX1,GX3,GX4,GX5,Xu1,Xu2,Xu3,XZ,XZ2} and others.
Following \cite{GX2}, the first author \cite{Xu3} introduced the following combinatorial $\alpha$-curvature for vertex scaling of PL metrics on triangulated surfaces.

\begin{definition}[\cite{Xu3}]\label{alpha_curvature definition}
Suppose $(S,V,\mathcal{T})$ is a triangulated surface with a discrete
PL metric $l$, $\alpha\in \mathbb{R}$ is a constant and $u: V\rightarrow \mathbb{R}$ is a discrete conformal factor for $l$.
The combinatorial $\alpha$-curvature at $i\in V$ is defined to be
\begin{equation*}
R_{\alpha, i}=\frac{K_i}{e^{\alpha u_i}},
\end{equation*}
where $K_i$ is the classical combinatorial curvature at $i\in V$ defined by (\ref{K_i}).
\end{definition}

\begin{remark}
If $\alpha=0$, the combinatorial $\alpha$-curvature $R_\alpha$ in Definition \ref{alpha_curvature definition}
is the classical discrete Gaussian curvature $K$.
Taking $g_i=e^{u_i}$ as a discrete analogue of the smooth Riemannian metric, we have
$R_{\alpha, i}(\lambda g_1,\cdots, \lambda g_V)=\lambda^{-\alpha}R_{\alpha, i}( g_1,\cdots, g_V)$ for any constant $\lambda>0$.
In the special case of $\alpha=1$, we have $R_{1, i}(\lambda g_1,\cdots, \lambda g_V)=\lambda^{-1}R_{1, i}(g_1,\cdots, g_V)$, which is parallelling to the transformation of smooth Gaussian curvature $K_{\lambda g}=\lambda^{-1} K_g$ with $g$ being the Riemannian metric.
\end{remark}

By the discrete Gauss-Bonnet formula (\ref{discrete Gauss-Bonnet formula}) for the classical discrete Gaussian curvature $K$,
the combinatorial $\alpha$-curvature $R_\alpha$ in Definition \ref{alpha_curvature definition} satisfies
the following discrete Gauss-Bonnet formula
$\sum_{i\in V}R_{\alpha,i}e^{\alpha u_i}=2\pi\chi(S).$
Therefore, if $\overline{R}\in \mathbb{R}^V$ is the combinatorial $\alpha$-curvature of some discrete PL metric discrete conformal to $l$
on $(S,V,\mathcal{T})$ with conformal factor $u$, then
\begin{equation}\label{prescribed alpha curvature formula}
\sum_{i\in V}\overline{R}_{i}e^{\alpha u_i}=2\pi\chi(S),
\end{equation}
which is a discrete analogue of the constraint equation $\int_S Ke^{2u}dV=2\pi\chi(S)$ in the smooth case \cite{KW1, Berger}.
Following Kazdan-Warner's arguments in \cite{KW1}, the constraint equation (\ref{prescribed alpha curvature formula}) imposes
the following sign conditions on $\overline{R}$ depending on $\chi(S)$:
\begin{description}\label{discrete KW condition}
  \item[(a)] $\chi(M)>0$: $\overline{R}$ is positive somewhere,
  \item[(b)] $\chi(M)=0$: $\overline{R}$ changes sign (unless $\overline{R}\equiv 0$),
  \item[(c)] $\chi(M)<0$: $\overline{R}$ is negative somewhere.
\end{description}
It is natural to ask the following discrete version of Kazdan-Warner's question for combinatorial $\alpha$-curvature.
\begin{question}[Discrete Kazdan-Warner Question]\label{discrete Kazdan-Warner question}
Suppose $(S, V)$ is a marked surface with a PL metric $d$.
  Are the sign conditions, depending on $\chi(S)$, sufficient conditions for a function $\overline{R}$ defined on $V$ to be
  the combinatorial $\alpha$-curvature of some polyhedral metric $d'$ discrete conformal to $d$?
\end{question}
We prove the following discrete Kazdan-Warner type theorem for Discrete Kazdan-Warner Question \ref{discrete Kazdan-Warner question}.
\begin{theorem}\label{main theorem 1}
Suppose $(S,V,d)$ is a marked surface with a PL metric $d$, $\alpha\in \mathbb{R}$ is a constant and $\overline{R}$ is a given function defined on $V$.
Then there exists a PL metric with combinatorial $\alpha$-curvature $\overline{R}$ in the discrete conformal class $\mathcal{D}(d)$ if one of the following conditions is satisfied :
\begin{description}
  \item[(1)] $\chi(S)>0,\ \alpha<0,\  \overline{R}>0$;
  \item[(2)] $\chi(S)<0,\ \alpha\neq0,\  \overline{R}\leq 0,\ \overline{R}\not\equiv 0$;
  \item[(3)] $\chi(S)=0,\ \alpha\neq0,\  \overline{R}\equiv0$;
  \item[(4)] $\alpha=0$, $\overline{R}\in (-\infty, 2\pi)$, $\sum_{i\in V}\overline{R}_{i}=2\pi \chi(S)$.
\end{description}
\end{theorem}

\begin{remark}
If $\overline{R}$ is a constant and $\alpha\overline{R}\leq0$, the existence of PL metric
with combinatorial $\alpha$-curvature $\overline{R}$ in the discrete conformal class $\mathcal{D}(d)$
has been proved in \cite{Xu3}.
In the case of $\alpha\overline{R}\leq0$, the uniqueness of PL metric with combinatorial $\alpha$-curvature $\overline{R}$ in the discrete conformal class $\mathcal{D}(d)$ has been proved in \cite{GLSW, Xu3}.
For the case $\alpha\overline{R}>0$, the uniqueness is unknown.
By the relationship of combinatorial $\alpha$-curvature and the classical discrete Gaussian curvature,
the cases \textbf{(3)} and \textbf{(4)} in Theorem \ref{main theorem 1}
are covered by Gu-Luo-Sun-Wu \cite{GLSW}.
Therefore, we just need to prove the cases \textbf{(1)} and \textbf{(2)} of Theorem \ref{main theorem 1}.
\end{remark}

The main tools for the proof of Theorem \ref{main theorem 1} are
Gu-Luo-Sun-Wu's discrete conformal theory for PL metrics on surfaces
\cite{GLSW} and variational principles with constraints.
The main ideas of the paper come from reading of Gu-Luo-Sun-Wu \cite{GLSW} and  Kou\v{r}imsk\'{a} \cite{Kourimska}.
A hyperbolic version of Theorem \ref{Euclidean discrete uniformization} has been proved by
Gu-Guo-Luo-Sun-Wu \cite{GGLSW}, which perfectly solves the classical prescribing discerete Gaussian curvature problem
in the hyperbolic background geometry. 
For prescribing combinatorial $\alpha$-curvature problem in the hyperbolic background geometry,
the authors \cite{XZ} recently obtained a hyperbolic version of Theorem \ref{main theorem 1}
using Luo's combinatorial Yamabe flow and Gu-Guo-Luo-Sun-Wu's discrete conformal theory for piecewise hyperbolic metrics
on surfaces \cite{GGLSW}.

\subsection{Organization of the paper}
The paper is organized as follows.
In Section \ref{section 2}, we recall the variational principle introduced by Luo \cite{Luo1} for vertex scaling and
Bobenko-Pinkall-Spingborn's development of Luo's variational principle,
and then recall the discrete conformal theory established by Gu-Luo-Sun-Wu \cite{GLSW}.
In Section \ref{section 3}, we translate Theorem \ref{main theorem 1} into an optimization problem with constraints.
In Section \ref{section 4}, we prove Theorem \ref{main theorem 1}.
\\
\\
\textbf{Acknowledgements}\\[8pt]
The first author thanks Professor Feng Luo, Dr. Tianqi Wu, Dr. Yanwen Luo,  Dr.  Wai Yeung Lam and
 Dr. Xiaoping Zhu for valuable communications.
The authors thank the organizers of
FRG Workshop on Geometric Methods for Analyzing Discrete Shapes for hosting the workshop that leads to this paper.
The research of the first author is supported by the Fundamental Research Funds for the Central Universities under
grant no. 2042020kf0199.

\section{Discrete conformality of PL metrics on surfaces}\label{section 2}
In this section, we recall some facts about discrete conformality of PL metrics that we need to use in the  proof
of Theorem \ref{main theorem 1}.

\subsection{The energy functions}

Suppose $\{ijk\}\in F_{\mathcal{T}}$ is a triangle and $l: E_\mathcal{T}\rightarrow \mathbb{R}_{>0}$
is a discrete PL metric on $(S, V, \mathcal{T})$.
Denote the inner angle in the triangle $\{ijk\}$ at the vertex $i$ as $\theta_i$.
Luo \cite{Luo1} proved the following result.

\begin{lemma}[\cite{Luo1}]\label{property on a triangle}
Suppose $\{ijk\}\in F_{\mathcal{T}}$ is a triangle and $l: E_\mathcal{T}\rightarrow \mathbb{R}_{>0}$
is a discrete PL metric on $(S, V, \mathcal{T})$.
\begin{description}
  \item[(1)] The admissible space of the discrete conformal factors for a triangle $\{ijk\}\in F_{\mathcal{T}}$
  \begin{equation*}
    \begin{aligned}
    \Omega_{ijk}=\{(u_i,u_j,u_k)\in \mathbb{R}^3| \widetilde{l}_{ij}, \widetilde{l}_{ik}, \widetilde{l}_{jk} \text{ satisfy the triangle inequality} \}
    \end{aligned}
  \end{equation*}
  is simply connected.
  \item[(2)] The Jacobian matrix $\frac{\partial (\theta_i, \theta_j, \theta_k)}{\partial (u_i, u_j, u_k)}$
  is symmetric and negative semi-definite with kernel $\{t(1,1,1)^T|t\in \mathbb{R}\}$ for any discrete conformal
  factor $(u_i, u_j, u_k)\in \Omega_{ijk}$.
\end{description}
\end{lemma}

Based on Lemma \ref{property on a triangle},
Luo \cite{Luo1} constructed the following energy function for a triangle $\{ijk\}\in F_{\mathcal{T}}$
\begin{equation}\label{energy function for a triangle}
    \begin{aligned}
    F_{ijk}(u_i,u_j,u_k)=-\int_{(0,0,0)}^{(u_i,u_j,u_k)}\theta_idu_i+\theta_jdu_j+\theta_kdu_k,
    \end{aligned}
\end{equation}
which is a well-defined smooth locally convex function defined on $\Omega_{ijk}$ with $\nabla F_{ijk}=(-\theta_i, -\theta_j, -\theta_k)$
and
\begin{equation}\label{property of F_ijk}
    \begin{aligned}
    F_{ijk}((u_i,u_j,u_k)+t(1,1,1))=F_{ijk}(u_i,u_j,u_k)-\pi t
    \end{aligned}
\end{equation}
for $t\in \mathbb{R}$.
Motivated by \cite{Verdiere}, Luo \cite{Luo1} further
introduced the following energy function
\begin{equation}\label{energy function for the surface}
    \begin{aligned}
   \mathbb{F}_\mathcal{T}(u)=\int_{0}^u\sum_{i\in V}K_idu_i=\sum_{\{ijk\}\in F_{\mathcal{T}}} F_{ijk}(u_i,u_j,u_k)+2\pi\sum_{i\in V}u_i,
    \end{aligned}
\end{equation}
which is a locally convex smooth function  of the discrete conformal factors
$u\in \cap_{\{ijk\}\in F_{\mathcal{T}}} \Omega_{ijk}$ with $\nabla \mathbb{F}_\mathcal{T}(u)=K$.
By (\ref{property of F_ijk}) and the Gauss-Bonnet formula,
the energy function $\mathbb{F}_\mathcal{T}(u)$ has  the following property
\begin{equation}\label{property of F}
\mathbb{F}_\mathcal{T}(u+c(1,\cdots,1))=\mathbb{F}_\mathcal{T}(u)+2c\pi\chi(S)
\end{equation}
for $c\in \mathbb{R}$  \cite{Luo1}.

In the following of the paper, we need to use the explicit expression of $\mathbb{F}_\mathcal{T}(u)$,
which was first constructed up to a constant by Bobenko-Pinkall-Spingborn \cite{BPS} as follows.
Set
\begin{equation*}
\Omega=\{(x,y,z)\in \mathbb{R}^3 |\ e^x+e^y>e^z, e^y+e^z>e^x, e^z+e^x>e^y\}.
\end{equation*}
For $(x,y,z)\in \Omega$, $e^x, e^y, e^z$ form the edge lengths of a Euclidean triangle.
Denote the inner angles facing $e^x, e^y, e^z$ as $\alpha, \beta, \gamma$ respectively and set
$$\mathbb{L}(x)=-\int_0^x \log|2\sin(t)|dt$$
to be Milnor's Lobachevsky function \cite{Milnor}.
Bobenko-Pinkall-Spingborn \cite{BPS} then defined the following function
\begin{equation*}
\begin{aligned}
f : \Omega &\rightarrow \mathbb{R}\\
(x,y,z)&\mapsto f(x,y,z)=\alpha x+\beta y+\gamma z+ \mathbb{L}(\alpha)+ \mathbb{L}(\beta)+ \mathbb{L}(\gamma).
\end{aligned}
\end{equation*}
Using $f$ as building blocks, Bobenko-Pinkall-Spingborn \cite{BPS} further constructed the following function
\begin{equation}\label{E T}
\begin{aligned}
\mathbb{E}_\mathcal{T}(u)
=\sum_{\{ijk\}\in F_\mathcal{T}}\left(2f(\frac{\widetilde{\lambda}_{ij}}{2}
,\frac{\widetilde{\lambda}_{jk}}{2},\frac{\widetilde{\lambda}_{ki}}{2})
-\frac{\pi}{2}(\widetilde{\lambda}_{ij}
+\widetilde{\lambda}_{jk}+\widetilde{\lambda}_{ki})\right)+2\pi\sum_{i\in V}u_i
\end{aligned}
\end{equation}
for $u\in \cap_{\{ijk\}\in F_{\mathcal{T}}} \Omega_{ijk}$, where $\widetilde{\lambda}_{ij}$
is the logarithm length of $\widetilde{l}_{ij}=l_{ij}e^{\frac{u_i+u_j}{2}}$.
Bobenko-Pinkall-Springborn \cite{BPS} proved that $\nabla \mathbb{E}_\mathcal{T}=K$,
which implies that $\mathbb{E}_\mathcal{T}(u)$ and $\mathbb{F}_\mathcal{T}(u)$ differs by some constant.

Bobenko-Pinkall-Springborn's construction of $\mathbb{E}_\mathcal{T}(u)$ in \cite{BPS}
is very elegant, but a little mysterious.
It seems that Bobenko-Pinkall-Springborn's construction
comes from their observation on the relationships of vertex scaling and 3-dimensional hyperbolic geometry
and can be taken as a consequence of the Schl\"{a}fli formula \cite{Rivin}.
Once one has the explicit form of the function $\mathbb{E}_\mathcal{T}(u)$,
it is easy to check that $\nabla \mathbb{E}_\mathcal{T}=K$.
However, the construction of $\mathbb{E}_\mathcal{T}(u)$ is not so easy.
As Luo's construction of the energy function $\mathbb{F}_\mathcal{T}(u)$ in \cite{Luo1} is relatively direct and easy,
a natural question of independent interest is whether one can derive an explicit form of $\mathbb{F}_\mathcal{T}(u)$ directly,
which differs from the explicit form of $\mathbb{E}_\mathcal{T}(u)$ by a constant.
In the following, we give such an argument,
the idea of which comes from Yuhao Xue from summer school 2017 at Tsinghua University.
In fact, we just need to derive an explicit form of the energy function $\mathbb{F}_{ijk}(u_i,u_j,u_k)$
for a triangle $\{ijk\}\in F_\mathcal{T}$.

\begin{proposition}
Suppose $\widetilde{l}$ is vertex scaling of $l$ on the triangle $\{ijk\}\in F_\mathcal{T}$
with a conformal factor $(u_i, u_j, u_k)\in \Omega_{ijk}$.
Denote the inner angle of the triangle with edge lengths $l_{ij}, l_{ik}, l_{jk}$
as $\overline{\Theta}_i, \overline{\Theta}_j, \overline{\Theta}_k$ respectively
and denote the inner angle of the triangle with edge lengths $\widetilde{l}_{ij}, \widetilde{l}_{ik}, \widetilde{l}_{jk}$
as $\Theta_i, \Theta_j, \Theta_k$ respectively.
Then the energy function $F_{ijk}(u_i,u_j,u_k)$ defined by (\ref{energy function for a triangle}) has the following explicit form
\begin{equation*}
\begin{aligned}
&F_{ijk}(u_i,u_j,u_k)\\
=&-(\Theta_iu_i+\Theta_ju_j+\Theta_ku_k)\\
  &+2\mathbb{L}(\Theta_i)+2\mathbb{L}(\Theta_j)+2\mathbb{L}(\Theta_k)
  +2\Theta_i\ln l_{jk}+2\Theta_j\ln l_{ik}+2\Theta_k\ln l_{ij}\\
  &-2\mathbb{L}(\overline{\Theta}_i)-2\mathbb{L}(\overline{\Theta}_j)-2\mathbb{L}(\overline{\Theta}_k)
  -2\overline{\Theta}_k \ln l_{ij}-2\overline{\Theta}_j\ln l_{ik}
  -2\overline{\Theta}_i\ln l_{jk}.
\end{aligned}
\end{equation*}
\end{proposition}
\proof
As $\Omega_{ijk}$ is simply connected by Lemma \ref{property on a triangle} and
$(0,0,0), (u_i, u_j, u_k)\in \Omega_{ijk}$, we can assume that $\gamma: [0,1]\rightarrow \Omega_{ijk}$
is a smooth path from $(0,0,0)$ to $(u_i, u_j, u_k)$.
Denote the corresponding quantities along the path $\gamma(t)$ as $u_i(t), \theta_i(t), l_{ij}(t)$ et al.
Then $u_i(1)=u_i, u_i(0)=0, \theta_i(1)=\Theta_i, \theta_i(0)=\overline{\Theta}_i$, $l_{ij}(1)=\widetilde{l}_{ij}$, $l_{ij}(0)=l_{ij}$.
By Definition \ref{luo's discrete conformal}, one can solve $u_i(t), u_j(t), u_k(t)$ as follows
\begin{equation}\label{u_i(t)}
\begin{aligned}
u_i(t)=\ln \frac{l_{ij}(t)l_{ik}(t)l_{jk}(0)}{l_{ij}(0)l_{ik}(0)l_{jk}(t)},
u_j(t)=\ln \frac{l_{ij}(t)l_{ik}(0)l_{jk}(t)}{l_{ij}(0)l_{ik}(t)l_{jk}(0)},
u_k(t)=\ln \frac{l_{ij}(0)l_{ik}(t)l_{jk}(t)}{l_{ij}(t)l_{ik}(0)l_{jk}(0)}.
\end{aligned}
\end{equation}
Suppose the circumcircle radius of the triangle with edge lengths $l_{ij}(t), l_{ik}(t), l_{jk}(t)$ is $R(t)$,
then
\begin{equation}\label{l_ij(t)}
\begin{aligned}
l_{ij}(t)=2R(t)\sin \theta_k(t), l_{ik}(t)=2R(t)\sin \theta_j(t), l_{jk}(t)=2R(t)\sin \theta_i(t).
\end{aligned}
\end{equation}
Submitting (\ref{l_ij(t)}) into  (\ref{u_i(t)}) gives
\begin{equation}\label{u_i(t) 2}
\begin{aligned}
u_i(t)=&3\ln R(t)+\ln \sin \theta_j(t)+\ln \sin \theta_k(t)-\ln \sin \theta_i(t)-\ln l_{ij}(0)-\ln l_{ik}(0)+\ln l_{jk}(0),\\
u_j(t)=&3\ln R(t)+\ln \sin \theta_i(t)+\ln \sin \theta_k(t)-\ln \sin \theta_j(t)-\ln l_{ij}(0)-\ln l_{jk}(0)+\ln l_{ik}(0),\\
u_k(t)=&3\ln R(t)+\ln \sin \theta_i(t)+\ln \sin \theta_j(t)-\ln \sin \theta_k(t)-\ln l_{ik}(0)-\ln l_{jk}(0)+\ln l_{ij}(0).
\end{aligned}
\end{equation}
By the definition of $F_{ijk}(u_i,u_j,u_k)$ in  (\ref{energy function for a triangle}), we have
\begin{equation}\label{F_ijk form proof 1}
\begin{aligned}
&F_{ijk}(u_i,u_j,u_k)\\
=&-\int_{0}^1\theta_i(t)du_i(t)+\theta_j(t)du_j(t)+\theta_k(t)du_k(t)\\
=&-[\theta_i(t)u_i(t)+\theta_j(t)u_j(t)+\theta_k(t)u_k(t)]|^1_0
    +\int_0^1u_i(t)d\theta_i(t)+u_j(t)d\theta_j(t)+u_k(t)d\theta_k(t)
\end{aligned}
\end{equation}
by integration by parts.
Submitting (\ref{u_i(t) 2}) into (\ref{F_ijk form proof 1}) gives
\begin{equation*}
\begin{aligned}
&F_{ijk}(u_i,u_j,u_k)\\
=&-[\theta_i(t)u_i(t)+\theta_j(t)u_j(t)+\theta_k(t)u_k(t)]|^1_0\\
  &-2\int_0^1[\theta_i'(t)\ln \sin \theta_i(t)+\theta_j'(t)\ln \sin \theta_j(t)+\theta_k'(t)\ln \sin \theta_k(t)]dt\\
  &+2\int_0^1[\ln l_{ij}(0)\theta_k'(t)+\ln l_{ik}(0)\theta_j'(t)+\ln l_{jk}(0)\theta_i'(t)]dt\\
=&-(\Theta_iu_i+\Theta_ju_j+\Theta_ku_k)\\
  &+2\mathbb{L}(\Theta_i)+2\mathbb{L}(\Theta_j)+2\mathbb{L}(\Theta_k)
  -2\mathbb{L}(\overline{\Theta}_i)-2\mathbb{L}(\overline{\Theta}_j)-2\mathbb{L}(\overline{\Theta}_k)\\
  &+2\ln l_{ij}(0)(\Theta_k-\overline{\Theta}_k)+2\ln l_{ik}(0)(\Theta_j-\overline{\Theta}_j)
  +2\ln l_{jk}(0)(\Theta_i-\overline{\Theta}_i),
\end{aligned}
\end{equation*}
where the identity $\theta_i'(t)+\theta_j'(t)+\theta_k'(t)\equiv0$ is used in the first equality.
\qed
\begin{remark}
  Set
\begin{equation*}
\begin{aligned}
h(u_i, u_j, u_j)=2f(\frac{\widetilde{\lambda}_{ij}}{2}
,\frac{\widetilde{\lambda}_{jk}}{2},\frac{\widetilde{\lambda}_{ki}}{2})
-\frac{\pi}{2}(\widetilde{\lambda}_{ij}
+\widetilde{\lambda}_{jk}+\widetilde{\lambda}_{ki}),
\end{aligned}
\end{equation*}
which is a building block of $\mathbb{E}_\mathcal{T}(u)$ in (\ref{E T}).
By direct calculations, we have
\begin{equation*}
\begin{aligned}
h(u_i, u_j, u_j)
=&-(\Theta_iu_i+\Theta_ju_j+\Theta_ku_k)\\
  &+2\mathbb{L}(\Theta_i)+2\mathbb{L}(\Theta_j)+2\mathbb{L}(\Theta_k)
  +2\Theta_i\ln l_{jk}+2\Theta_j\ln l_{ik}+2\Theta_k\ln l_{ij}\\
  &-\pi(\ln l_{ij}+\ln l_{ik}+\ln l_{jk}),
\end{aligned}
\end{equation*}
which implies that $F_{ijk}(u_i,u_j,u_k)$ differ from $h(u_i, u_j, u_j)$ by a constant depending on the background discrete PL metric $l$.
As a direct consequence, Luo's energy function $\mathbb{F}_\mathcal{T}(u)$ in (\ref{energy function for the surface})
differs from Bobenko-Pinkall-Springborn's energy function $\mathbb{E}_\mathcal{T}(u)$ in (\ref{E T})
by a constant depending on the background discrete PL  metric $l$.
\end{remark}

\begin{remark}
There are two approaches to extend the energy function
$\mathbb{F}_\mathcal{T}(u)$ or $\mathbb{E}_\mathcal{T}(u)$ to be a convex function
defined on $\mathbb{R}^V$ introduced by Bobenko-Pinkall-Springborn \cite{BPS} and Gu-Luo-Sun-Wu \cite{GLSW} independently.
In \cite{BPS}, the function $\mathbb{E}_\mathcal{T}(u)$ is extended to be
a $C^1$ smooth convex function defined on $\mathbb{R}^V$
by exploiting the fact that the function $f$ defined on $\Omega$ could be
extended to a $C^1$ smooth convex function on $\mathbb{R}^3$ by extending the inner angles of a triangle by constants.
Using Luo' development \cite{Luo3} of Bobenko-Pinkall-Spingborn's extension \cite{BPS},
Ge-Jiang \cite{Ge-Jiang} similarly extended
the potential function $\mathbb{F}_\mathcal{T}(u)$ defined by (\ref{energy function for the surface})
to be a $C^1$-smooth convex function defined on $\mathbb{R}^V$.
Luo's development \cite{Luo3} of Bobenko-Pinkall-Spingborn's extension \cite{BPS} has recently been further developed to
handle other cases. See, for instance, \cite{Xu2,Xu4,Xu5,XZ0} and others.
In \cite{GLSW}, the function $\mathbb{F}_\mathcal{T}(u)$ is extended
to be a $C^2$-smooth convex function defined on $\mathbb{R}^V$
 by changing the triangulation of the marked surface under the Delaunay condition,
 based on which Gu-Luo-Sun-Wu proved the Discrete Uniformization Theorem \ref{Euclidean discrete uniformization}
 for PL metrics on closed surfaces.
The first approach could not ensure the triangles being non-degenerate, while the second approach could.
\end{remark}

\subsection{Gu-Luo-Sun-Wu's discrete conformal theory}

Based on Penner's decorated Teichim\"{u}ller space theory \cite{Penner},
Gu-Luo-Sun-Wu \cite{GLSW} established the discrete conformal theory for PL metrics on compact surfaces and
proved Theorem \ref{Euclidean discrete uniformization}.
In the following, we recall some results in \cite{GLSW} that we need.
We will not  formally involve too much notions and results on decorated Teichim\"{u}ller space.
For more details of these results, please refer to Gu-Luo-Sun-Wu's original work \cite{GLSW}.

\begin{theorem}[\cite{GLSW} Corollary 4.7]\label{D(d) parameterized by Rv}
Suppose $d$ is a PL metric on the marked surface $(S, V)$.
Then there exists a $C^1$ diffeomorphism $\phi: \mathcal{D}(d)\rightarrow \mathbb{R}^V$.
\end{theorem}

By Theorem \ref{D(d) parameterized by Rv}, the discrete conformal class $\mathcal{D}(d)$
is parameterized by $\mathbb{R}^V$. For simplicity, for any $d'\in \mathcal{D}(d)$,
we denote it by $d(u)$ for some $u\in \mathbb{R}^V$.
Suppose $\mathcal{T}$ is a triangulation of the marked surface $(S,V)$.
Set
\begin{equation*}
\mathcal{A}_\mathcal{T}=\{u\in R^V |\ \mathcal{T}\ \text{is isotopic to a Delaunay triangulation of}\ (S,V,d(u))\}.
\end{equation*}
Based on Akiyoshi's finiteness theorem in \cite{Akiyoshi} (see also Appendix in \cite{GLSW} for a new proof),
Gu-Luo-Sun-Wu \cite{GLSW} proved the following result.

\begin{theorem}[\cite{GLSW} Lemma 5.1]\label{finite cell decomposition}
Let
$$J=\{\mathcal{T}| \mathcal{A}_\mathcal{T}\ \text{has nonempty interior in}\ \mathbb{R}^V\}.$$
Then $J$ is a finite set, $\mathbb{R}^V=\cup_{\mathcal{T}_i\in J}\mathcal{A}_{\mathcal{T}_i}$
and $\mathcal{A}_{\mathcal{T}_i}$ is real analytically diffeomorphic to a closed convex polytope in $\mathbb{R}^V$.
\end{theorem}

By Theorem \ref{D(d) parameterized by Rv}, Gu-Luo-Sun-Wu \cite{GLSW} introduced the following globally defined  $C^1$ smooth
combinatorial curvature function
\begin{equation*}
\begin{aligned}
\mathbf{F}: \mathbb{R}^V&\rightarrow (-\infty, 2\pi)^V\\
u&\mapsto K(d(u)).
\end{aligned}
\end{equation*}
Combining with Lemma \ref{property on a triangle}, Gu-Luo-Sun-Wu \cite{GLSW} further constructed
a globally defined energy function with $\mathbf{F}$ as gradient,
which plays an important role in the proof of Theorem \ref{Euclidean discrete uniformization}.

\begin{theorem}[\cite{GLSW} Proposition 5.2]\label{thm global energy function}
  There exists a $C^2$-smooth convex function
  \begin{equation}\label{globally defined energy function W}
\begin{aligned}
\mathbb{E}: \mathbb{R}^V&\rightarrow \mathbb{R}\\
u&\mapsto \int_0^u\sum_{i\in V}\mathbf{F}_idu_i
\end{aligned}
\end{equation}
so that its gradient $\nabla \mathbb{E}=\mathbf{F}$ and the restriction of $\mathbb{E}$
to the hyperplane $\{u\in \mathbb{R}^V|\sum_{i\in V}u_i=0\}$ is strictly convex.
\end{theorem}

\begin{remark}\label{relation of W and F_T}
By the construction in Theorem \ref{thm global energy function}, for $\mathcal{T}\in J$, the restriction
$\mathbb{E}|_{\mathcal{A}_\mathcal{T}}$ differs from $\mathbb{F}_\mathcal{T}$ in (\ref{energy function for the surface})
by a constant.
As a consequence of (\ref{property of F}), $\mathbb{E}$ has the following property
\begin{equation}\label{property of E globally}
\mathbb{E}(u+c(1,\cdots,1))=\mathbb{E}(u)+2c\pi\chi(S)
\end{equation}
for any $c\in \mathbb{R}$.
\end{remark}

\section{Variational principles with constraints}\label{section 3}
In this section, we translate Theorem \ref{main theorem 1} into an optimization problem with inequality constraints by variational principles, which involve the function $\mathbb{E}$ defined in (\ref{globally defined energy function W}).

Suppose $(S,V,d)$ is a marked surface with a PL metric $d$,
$\alpha\in \mathbb{R}$ is a non-zero constant and $\overline{R}$ is a given function defined on $V$.
Set
\begin{equation}\label{A}
\mathcal{A}=\{u\in \mathbb{R}^V|0>\sum_{i\in V}\overline{R}_i e^{\alpha u_i}\geq 2\pi\chi(S),\  \overline{R}\leq0,\ \overline{R}\not\equiv0\},
\end{equation}
\begin{equation}\label{B}
\mathcal{B}=\{u\in \mathbb{R}^V|0<\sum_{i\in V}\overline{R}_i e^{\alpha u_i}\leq 2\pi\chi(S),\ \overline{R}>0\},
\end{equation}
\begin{equation}\label{C}
\mathcal{C}=\{u\in \mathbb{R}^V|\sum_{i\in V} \overline{R}_i e^{\alpha u_i}\leq 2\pi\chi(S)<0,\ \overline{R}\leq0,\ \overline{R}\not\equiv0\}.
\end{equation}

\begin{proposition}\label{ABC proposition}
The sets $\mathcal{A},\ \mathcal{B}$ and $\mathcal{C}$ are unbounded closed subsets of $\mathbb{R}^V$.
\end{proposition}
\proof
It is obvious that the sets $\mathcal{A},\ \mathcal{B}$ and $\mathcal{C}$ are closed subsets of $\mathbb{R}^V$.
For the set $\mathcal{A}$, by direct calculations,
\begin{equation*}
\sum_{i\in V} \overline{R}_i e^{\alpha (u_i+c)}=e^{\alpha c}\sum_{i\in V} \overline{R}_i e^{\alpha u_i}\geq 2\pi\chi(S)
\end{equation*}
is equivalent to
\begin{equation*}
c\geq\frac{1}{\alpha}\log\frac{2\pi\chi(S)} {\sum_{i\in V} \overline{R}_i e^{\alpha u_i}}
\end{equation*}
under the condition $\alpha<0$;
\begin{equation*}
\sum_{i\in V} \overline{R}_i e^{\alpha (u_i+c)}=e^{\alpha c}\sum_{i\in V} \overline{R}_i e^{\alpha u_i}\geq 2\pi\chi(S)
\end{equation*}
is equivalent to
\begin{equation*}
c\leq\frac{1}{\alpha}\log\frac{2\pi\chi(S)} {\sum_{i\in V} \overline{R}_i e^{\alpha u_i}}
\end{equation*}
under the condition $\alpha>0$.
Therefore, the set $\mathcal{A}$ is unbounded. Similarly, the sets $\mathcal{B}$ and $\mathcal{C}$ are unbounded.
\qed

According to Proposition \ref{ABC proposition}, we have following result.
\begin{lemma}\label{minimum lies at the boundary}
Suppose $(S,V,d)$ is a marked surface with a PL metric $d$, $\alpha\in \mathbb{R}$ is a constant and $\overline{R}$ is a given function defined on $V$. If one of the following three conditions is satisfied
\begin{description}
  \item[(1)] $\alpha>0$ and the energy function $\mathbb{E}$ attains a minimum in the set $\mathcal{A}$,
  \item[(2)] $\alpha<0$ and the energy function $\mathbb{E}$ attains a minimum in the set $\mathcal{B}$,
  \item[(3)] $\alpha<0$ and the energy function $\mathbb{E}$ attains a minimum in the set $\mathcal{C}$,
\end{description}
then the minimum value point of  $\mathbb{E}$ lies in the set $\{u\in \mathbb{R}^V|\sum_{i\in V} \overline{R}_i e^{\alpha u_i}=2\pi\chi(S)\}$.
\end{lemma}
\proof
Let $\alpha>0$ and suppose the function $\mathbb{E}$ attains a minimum at $u\in \mathcal{A}$.
The definition of $\mathcal{A}$ in (\ref{A}) implies $\chi(S)<0$.
Taking $c_0=\frac{1}{\alpha}\log\frac{2\pi\chi(S)} {\sum_{i\in V} \overline{R}_i e^{\alpha u_i}}$, then $c_0\geq0$.
By the proof of Proposition \ref{ABC proposition}, $u+c_0\mathbb{I}\in \mathcal{A}$.
Therefore, by the additive property of the function $\mathbb{E}$ in (\ref{property of E globally}),
we have
\begin{equation*}
\mathbb{E}(u)\leq \mathbb{E}(u+c_0\mathbb{I})=\mathbb{E}(u)+2\pi c_0\chi(S),
\end{equation*}
which implies $c_0\leq0$ by $\chi(S)<0$.
Hence $c_0=0$ and $\sum_{i\in V} \overline{R}_i e^{\alpha u_i}=2\pi\chi(S)$.
This proves the case of $\mathbf{(1)}$.
The proofs for the cases $\mathbf{(2)}$ and $\mathbf{(3)}$ are similar, we omit the details here.
\qed

By Lemma \ref{minimum lies at the boundary}, we translate Theorem \ref{main theorem 1} into the following theorem, which is a non-convex optimization problem with inequality constraints.

\begin{theorem}\label{inequality constraints theorem}
Suppose $(S,V,d)$ is a marked surface with a PL metric $d$ and $\chi(S)\neq0$, $\alpha\in \mathbb{R}$ is a non-zero constant and $\overline{R}$ is a given function defined on $V$.
\begin{description}
  \item[(1)] If $\alpha>0$ and the energy function $\mathbb{E}$ attains a minimum in $\mathcal{A}$, then there exists a PL metric in the conformal class $\mathcal{D}(d)$ with combinatorial $\alpha$-curvature $\overline{R}\leq0$ and $\overline{R}\not\equiv0$;
  \item[(2)]If $\alpha<0$ and the energy function $\mathbb{E}$ attains a minimum in $\mathcal{B}$, then there exists a PL metric in the conformal class $\mathcal{D}(d)$ with combinatorial $\alpha$-curvature $\overline{R}>0$;
  \item[(3)] If $\alpha<0$ and the energy function $\mathbb{E}$ attains a minimum in $\mathcal{C}$, then there exists a PL metric in the conformal class $\mathcal{D}(d)$ with combinatorial $\alpha$-curvature $\overline{R}\leq0$ and $\overline{R}\not\equiv0$.
\end{description}

\end{theorem}
\proof
Lemma \ref{minimum lies at the boundary} shows that if $u\in \mathbb{R}^V$ is a minimum of the energy
function $\mathbb{E}$ defined on one of these sets, then $\sum_{i\in V} \overline{R}_i e^{\alpha u_i}= 2\pi\chi(S)$.
The conclusion follows from the following claim.\\
\textbf{Claim :} Up to scaling, the PL-metrics with prescribed combinatorial $\alpha$-curvature in the discrete conformal class $\mathcal{D}(d)$ are in one-to-one correspondence with the critical points of the function $\mathbb{E}$ under the constraint $\sum_{i\in V} \overline{R}_i e^{\alpha u_i}=2\pi\chi(S)$.

We use the method of Lagrange multipliers to prove this claim.
Set
\begin{equation*}
H(u,\lambda)=\mathbb{E}(u)+\lambda \left(\sum_{i\in V} \overline{R}_i e^{\alpha u_i}-2\pi\chi(S)\right),
\end{equation*}
where $\lambda\in \mathbb{R}$ is a Lagrange multiplier.
If $u$ is a critical point of the function $\mathbb{E}$ under the constraint $\sum_{i\in V} \overline{R}_i e^{\alpha u_i}=2\pi\chi(S)$,
then
\begin{equation*}
0=\frac{\partial H(u,\lambda)}{\partial u_i}=K_i+\lambda\alpha\overline{R}_i e^{\alpha u_i},
\end{equation*}
which implies
\begin{equation*}
R_{\alpha, i}=\frac{K_i}{e^{\alpha u_i}}=-\lambda\alpha\overline{R}_i.
\end{equation*}
By the discrete Gauss-Bonnet formula (\ref{discrete Gauss-Bonnet formula}), the Lagrange multiplier $\lambda$ satisfies
\begin{equation*}
\lambda=-\frac{2\pi \chi(S)}{\alpha\sum_{i\in V} \overline{R}_i e^{\alpha u_i}}=-\frac{1}{\alpha}
\end{equation*}
under the constraint
$\sum_{i\in V} \overline{R}_i e^{\alpha u_i}=2\pi\chi(S)$,
which implies the combinatorial $\alpha$-curvature
\begin{equation*}
R_{\alpha, i}=-\lambda\alpha\overline{R}_i=\frac{2\pi \chi(S)}{\sum_{i\in V} \overline{R}_i e^{\alpha u_i}}\overline{R}_i=\overline{R}_i
\end{equation*}
under the constraint
$\sum_{i\in V} \overline{R}_i e^{\alpha u_i}=2\pi\chi(S)$.
\qed

\section{Proof of Theorem \ref{main theorem 1}}\label{section 4}

In this section, we give a proof for Theorem \ref{main theorem 1}.
By Theorem \ref{inequality constraints theorem},
we just need to prove that the energy function $\mathbb{E}$ attains the minimum
in the sets $\mathcal{A},\ \mathcal{B}$ and $\mathcal{C}$.
Recall the following classical result from calculus.

\begin{theorem}\label{calculus theorem}
Let $A\subseteq \mathbb{R}^m$ be a closed set and $f: A\rightarrow \mathbb{R}$ be a continuous function. If every unbounded sequence $\{u_n\}_{n\in \mathbb{N}}$ in $A$ has a subsequence $\{x_{n_k}\}_{k\in \mathbb{N}}$ such that
$\lim_{k\rightarrow +\infty} f(x_{n_k})=+\infty$,
then $f$ attains a minimum in $A$.
\end{theorem}

One can refer to \cite{Kourimska Thesis} (Section 4.1) for a proof of Theorem \ref{calculus theorem}.
The majority of the conditions in Theorem \ref{calculus theorem} are satisfied,
including the sets $\mathcal{A},\ \mathcal{B}$ and $\mathcal{C}$ are closed subsets of $\mathbb{R}^V$ by Proposition \ref{ABC proposition} and the energy function $\mathbb{E}$ is continuous by Theorem \ref{thm global energy function}.
To prove Theorem \ref{main theorem 1},
we just need to prove the following theorem by Theorem \ref{inequality constraints theorem} and Theorem \ref{calculus theorem}.

\begin{theorem}\label{main theorem 2}
Suppose $(S,V,d)$ is a marked surface with a PL metric $d$, $\alpha\in \mathbb{R}$ is a constant and $\overline{R}$ is a given function defined on $V$. If one of the following three conditions is satisfied,
\begin{description}
  \item[(1)] $\alpha>0$ and $\{u_n\}_{n\in \mathbb{N}}$ is an unbounded sequence in $\mathcal{A}$,
  \item[(2)] $\alpha<0$ and $\{u_n\}_{n\in \mathbb{N}}$ is an unbounded sequence in $\mathcal{B}$,
  \item[(3)] $\alpha<0$ and $\{u_n\}_{n\in \mathbb{N}}$ is an unbounded sequence in $\mathcal{C}$,
\end{description}
then there exist a subsequence $\{u_{n_k}\}_{k\in \mathbb{N}}$ of $\{u_n\}_{n\in \mathbb{N}}$ such that
$\lim_{k\rightarrow +\infty} \mathbb{E}(u_{n_k})=+\infty$.
\end{theorem}

Let $\{u_n\}_{n\in \mathbb{N}}$ be an unbounded sequence in $\mathbb{R}^V$,
denote its coordinate sequence at $j\in V$ by $\{u_{j,n}\}_{n\in \mathbb{N}}$.
Following \cite{Kourimska}, we use the following convention.

\begin{convention}\label{convention}
The sequence $\{u_n\}_{n\in \mathbb{N}}$ possesses the following properties:
\begin{description}
  \item[(1)] It lies in one cell $\mathcal{A}_\mathcal{T}$ of $\mathbb{R}^V$ for some $\mathcal{T}\in J$ given by
  Theorem \ref{finite cell decomposition};
  \item[(2)] There exists a vertex $i^*\in V$ such that $u_{i^*,n}\leq u_{j,n}$ for all $j\in V$ and $n\in \mathbb{N}$;
  \item[(3)] Each coordinate sequence $\{u_{j,n}\}_{n\in \mathbb{N}}$ either converge, diverge properly to $+\infty$, or diverges properly to $-\infty$;
  \item[(4)] For all $j\in V$, the sequence $\{u_{j,n}-u_{i^*,n}\}_{n\in \mathbb{N}}$ either converge or diverge properly to $+\infty$.
\end{description}
\end{convention}

By Theorem \ref{finite cell decomposition}, it is obvious that every sequence in $\mathbb{R}^V$ possesses a subsequence satisfying these properties in Convention \ref{convention}, hence the sequence in Convention \ref{convention} could be chosen without loss of generality.
To prove Theorem \ref{main theorem 2}, we further need the following two results obtained by Kou\v{r}imsk\'{a} \cite{Kourimska}.

\begin{lemma}[\cite{Kourimska} Corollary 5.6]\label{triangle sequence converge lemma}
In every triangle $\{ijk\}\in F_\mathcal{T}$, at least two of the three sequences $(u_{i,n}-u_{i^*,n})_{n\in \mathbb{N}}$, $(u_{j,n}-u_{i^*,n})_{n\in \mathbb{N}}$ and $(u_{k,n}-u_{i^*,n})_{n\in \mathbb{N}}$ converge.
\end{lemma}

\begin{lemma}[\cite{Kourimska} Lemma 5.11]\label{E decomposition lemma}
There exists a convergent sequence $\{C_n\}_{n\in \mathbb{N}}$ such that the energy function $\mathbb{E}$ satisfies
\begin{equation*}
\mathbb{E}(u_n)=C_n+2\pi\left(u_{i^*,n}\chi(S)+\sum_{j\in V}(u_{j,n}-u_{i^*,n})\right).
\end{equation*}
\end{lemma}
The proof of Lemma \ref{E decomposition lemma} is based on an interesting analysis of the explicit form of the energy function $\mathbb{F}_\mathcal{T}$ or $\mathbb{E}_\mathcal{T}$. Readers can refer to \cite{Kourimska,Kourimska Thesis} for the proof.\\
\bigskip
\\
\noindent\textbf{Proof of Theorem \ref{main theorem 2}:}
Assume that $\{u_n\}_{n\in \mathbb{N}}$ is an unbounded sequence satisfying Convention \ref{convention}.
We just need to prove that $\lim_{n\rightarrow +\infty} \mathbb{E}(u_n)=+\infty$.
\begin{description}
  \item[(1)]
Let $\alpha>0$ and $\{u_n\}_{n\in \mathbb{N}}$ is an unbounded sequence in $\mathcal{A}$.
The definition of $\mathcal{A}$ in (\ref{A}) implies $\chi(S)<0$, $\overline{R}\leq0$ and $\overline{R}\not\equiv0$.
Since the sequence $\{u_n\}_{n\in \mathbb{N}}$ lies in $\mathcal{A}$, we have
\begin{equation}\label{(1)}
0>\sum_{j\in V} \overline{R}_j e^{\alpha (u_{j,n}-u_{i^*,n})}=e^{-\alpha u_{i^*,n}}\cdot\sum_{j\in V} \overline{R}_j e^{\alpha u_{j,n}}\geq 2\pi\chi(S) e^{-\alpha u_{i^*,n}}.
\end{equation}
Note that $\left(\sum_{j\in V}(u_{j,n}-u_{i^*,n})\right)_{n\in \mathbb{N}}$ converges to a finite number or
diverges properly  to $+\infty$ by Convention \ref{convention}.

If $\left(\sum_{j\in V}(u_{j,n}-u_{i^*,n})\right)_{n\in \mathbb{N}}$ converges to a finite number,
then the sequence $(u_{j,n}-u_{i^*,n})_{n\in \mathbb{N}}$ converges for all $j\in V$ by Convention \ref{convention},
which implies
$\sum_{j\in V} \overline{R}_j e^{\alpha (u_{j,n}-u_{i^*,n})}$ converges to a finite negative number by
$\overline{R}\leq0$ and $\overline{R}\not\equiv0$.
Therefore, $\{u_{i^*,n}\}_{n\in \mathbb{N}}$ is bounded from above by (\ref{(1)}), $\alpha>0$ and $\chi(S)<0$,
which implies $\{u_{i^*,n}\}_{n\in \mathbb{N}}$ converges to a finite number or diverges to $-\infty$.
If $\{u_{i^*,n}\}_{n\in \mathbb{N}}$ converges to a finite number,
then $\{u_{j,n}\}_{n\in \mathbb{N}}$ are bounded for all $j\in V$ by
$(u_{j,n}-u_{i^*,n})_{n\in \mathbb{N}}$ converges for all $j\in V$, which implies $\{u_n\}_{n\in \mathbb{N}}$ is bounded.
This contradicts the assumption that $\{u_n\}_{n\in \mathbb{N}}$ is unbounded.
Therefore, the sequence $\{u_{i^*,n}\}_{n\in \mathbb{N}}$ diverges properly to $-\infty$.
Combining this with $\chi(S)<0$ and Lemma \ref{E decomposition lemma},
we have $\lim_{n\rightarrow +\infty} \mathbb{E}(u_n)=+\infty$.

If $\left(\sum_{j\in V}(u_{j,n}-u_{i^*,n})\right)_{n\in \mathbb{N}}$ diverges properly to $+\infty$,
then there exists at least one vertex $j\in V$ such that the sequence $(u_{j,n}-u_{i^*,n})_{n\in \mathbb{N}}$ diverges properly to $+\infty$.
By Lemma \ref{triangle sequence converge lemma},
there exists at least one vertex $k\in V$ such that
the sequence $(u_{k,n}-u_{i^*,n})_{n\in \mathbb{N}}$ converges (for example $(u_{i^*,n}-u_{i^*,n})_{n\in \mathbb{N}}$).
Therefore, $e^{\alpha (u_{j,n}-u_{i^*,n})}$ converges to a finite positive number or diverges properly to $+\infty$ and
for at least one vertex $j\in V$ the term $e^{\alpha (u_{j,n}-u_{i^*,n})}$ converges to a finite positive number.
Combining $\overline{R}\leq0$ and $\overline{R}\not\equiv0$,
$\sum_{j\in V} \overline{R}_j e^{\alpha (u_{j,n}-u_{i^*,n})}$
tends to either a finite negative number or $-\infty$.
If $\sum_{j\in V} \overline{R}_j e^{\alpha (u_{j,n}-u_{i^*,n})}$ tends to a finite negative number,
then $u_{i^*,n}$ is bounded from above by (\ref{(1)}), which implies that $u_{i^*,n}\chi(S)$ is bounded from below.
Combining with the assumption that
$\left(\sum_{j\in V}(u_{j,n}-u_{i^*,n})\right)_{n\in \mathbb{N}}$ diverges properly to $+\infty$, we have
$\lim_{n\rightarrow +\infty} \mathbb{E}(u_n)=+\infty$ by Lemma \ref{E decomposition lemma}.
If $\sum_{j\in V} \overline{R}_j e^{\alpha (u_{j,n}-u_{i^*,n})}$ tends to $-\infty$,
we have $\{u_{i^*,n}\}_{n\in \mathbb{N}}$ diverges properly to $-\infty$ by (\ref{(1)}), which
implies $\lim_{n\rightarrow +\infty} \mathbb{E}(u_n)=+\infty$ by $\chi(S)<0$ and Lemma \ref{E decomposition lemma}.

\item[(2)] Let $\alpha<0$ and $\{u_n\}_{n\in \mathbb{N}}$ is an unbounded sequence in $\mathcal{B}$. The definition of $\mathcal{B}$ in (\ref{B}) implies $\chi(S)>0$ and $\overline{R}>0$.
Since the sequence $\{u_n\}_{n\in \mathbb{N}}$ lies in $\mathcal{B}$, we have
\begin{equation}\label{(2)}
0<\sum_{j\in V} \overline{R}_j e^{\alpha (u_{j,n}-u_{i^*,n})}=e^{-\alpha u_{i^*,n}}\cdot\sum_{j\in V} \overline{R}_j e^{\alpha u_{j,n}}\leq 2\pi\chi(S) e^{-\alpha u_{i^*,n}}.
\end{equation}
If $\left(\sum_{j\in V}(u_{j,n}-u_{i^*,n})\right)_{n\in \mathbb{N}}$ converges, then the sequence $(u_{j,n}-u_{i^*,n})_{n\in \mathbb{N}}$ converges for all $j\in V$ by Convention \ref{convention}.
Note that $\alpha<0$, $\chi(S)>0$ and $\overline{R}>0$,
we have $\{u_{i^*,n}\}_{n\in \mathbb{N}}$ is bounded from below by (\ref{(2)}),
which implies $\{u_{i^*,n}\}_{n\in \mathbb{N}}$ converges to a finite number or
diverges properly to $+\infty$ by Convention \ref{convention}.
Combining this with $\{u_n\}_{n\in \mathbb{N}}$ is unbounded
and $(u_{j,n}-u_{i^*,n})_{n\in \mathbb{N}}$ converges for all $j\in V$,
we have the sequence $\{u_{i^*,n}\}_{n\in \mathbb{N}}$ diverges properly to $+\infty$.
As a result, we have $\lim_{n\rightarrow +\infty} \mathbb{E}(u_n)=+\infty$
by Lemma \ref{E decomposition lemma} and $\chi(S)>0$.

If the sequence $\left(\sum_{j\in V}(u_{j,n}-u_{i^*,n})\right)_{n\in \mathbb{N}}$ diverges properly to $+\infty$,
then there exists at least one vertex $j\in V$ such that the sequence $(u_{j,n}-u_{i^*,n})_{n\in \mathbb{N}}$ diverges properly to $+\infty$.
By Lemma \ref{triangle sequence converge lemma},
there exists at least one vertex $k\in V$ such that
the sequence $(u_{k,n}-u_{i^*,n})_{n\in \mathbb{N}}$ converges (for example $(u_{i^*,n}-u_{i^*,n})_{n\in \mathbb{N}}$).
Therefore, $\sum_{j\in V} \overline{R}_j e^{\alpha (u_{j,n}-u_{i^*,n})}$
has a positive lower bound by $\overline{R}>0$ and $\alpha<0$, which implies
$2\pi\chi(S) e^{-\alpha u_{i^*,n}}$ has a positive lower bound by (\ref{(2)}).
Therefore, $\{u_{i^*,n}\}_{n\in \mathbb{N}}$ is bounded from below by $\alpha<0$ and $\chi(S)>0$ and hence $u_{i^*,n}\chi(S)$ is bounded from below.
Combining this with $\left(\sum_{j\in V}(u_{j,n}-u_{i^*,n})\right)_{n\in \mathbb{N}}$ diverges properly to $+\infty$,
we have $\lim_{n\rightarrow +\infty} \mathbb{E}(u_n)=+\infty$ by Lemma \ref{E decomposition lemma}.

\item[(3)]
Let $\alpha<0$ and $\{u_n\}_{n\in \mathbb{N}}$ is an unbounded sequence in $\mathcal{C}$.
The definition of $\mathcal{C}$ in ($\ref{C}$) implies $\chi(S)<0$, $\overline{R}\leq0$ and $\overline{R}\not\equiv0$.
Since the sequence $\{u_n\}_{n\in \mathbb{N}}$ lies in $\mathcal{C}$, we have
\begin{equation}\label{(3)}
\sum_{j\in V} \overline{R}_j e^{\alpha (u_{j,n}-u_{i^*,n})}=e^{-\alpha u_{i^*,n}}\cdot\sum_{j\in V} \overline{R}_j e^{\alpha u_{j,n}}\leq 2\pi\chi(S)e^{-\alpha u_{i^*,n}}<0.
\end{equation}
If $\left(\sum_{j\in V}(u_{j,n}-u_{i^*,n})\right)_{n\in \mathbb{N}}$ converges, then the sequence $(u_{j,n}-u_{i^*,n})_{n\in \mathbb{N}}$ converges for all $j\in V$ by Convention \ref{convention},
which implies that
$\sum_{j\in V} \overline{R}_j e^{\alpha (u_{j,n}-u_{i^*,n})}$ converges to a finite negative number
by $\overline{R}\leq0$ and $\overline{R}\not\equiv0$.
Combining this with $\alpha<0$ and $\chi(S)<0$,
we have $\{u_{i^*,n}\}_{n\in \mathbb{N}}$ is bounded from above by (\ref{(3)}).
As $\{u_n\}_{n\in \mathbb{N}}$ is unbounded, then similar to the arguments above,
the sequence $\{u_{i^*,n}\}_{n\in \mathbb{N}}$ diverges properly to $-\infty$.
Combining this with $\chi(S)<0$ and Lemma \ref{E decomposition lemma},
we have $\lim_{n\rightarrow +\infty} \mathbb{E}(u_n)=+\infty$.

If $\left(\sum_{j\in V}(u_{j,n}-u_{i^*,n})\right)_{n\in \mathbb{N}}$ diverges properly to $+\infty$,
then there exists at least one vertex $j\in V$
such that the sequence $(u_{j,n}-u_{i^*,n})_{n\in \mathbb{N}}$ diverges properly to $+\infty$.
By Lemma \ref{triangle sequence converge lemma},
there exists at least one vertex $k\in V$ such that
the sequence $(u_{k,n}-u_{i^*,n})_{n\in \mathbb{N}}$ converges (for example $(u_{i^*,n}-u_{i^*,n})_{n\in \mathbb{N}}$).
Therefore, $\sum_{j\in V} \overline{R}_j e^{\alpha (u_{j,n}-u_{i^*,n})}$ either tends to zero or to a finite negative number
by $\alpha<0$, $\overline{R}\leq0$ and $\overline{R}\not\equiv0$.
If $\sum_{j\in V} \overline{R}_j e^{\alpha (u_{j,n}-u_{i^*,n})}$ tends to zero, then $2\pi\chi(S) e^{-\alpha u_{i^*,n}}$ tends to zero by (\ref{(3)}), which implies $\{u_{i^*,n}\}_{n\in \mathbb{N}}$ diverges properly to $-\infty$ by $\alpha<0$ and $\chi(S)<0$.
Combining this with $\chi(S)<0$ and Lemma \ref{E decomposition lemma},
we have $\lim_{n\rightarrow +\infty} \mathbb{E}(u_n)=+\infty$.
If $\sum_{j\in V} \overline{R}_j e^{\alpha (u_{j,n}-u_{i^*,n})}$ tends to a finite negative number,
we have $2\pi\chi(S) e^{-\alpha u_{i^*,n}}$ has a negative lower bound by (\ref{(3)}),
which implies $u_{i^*,n}$ is bounded from above by $\alpha<0$ and $\chi(S)<0$.
Combining this with $\chi(S)<0$ and $\left(\sum_{j\in V}(u_{j,n}-u_{i^*,n})\right)_{n\in \mathbb{N}}$
diverges properly to $+\infty$,
we have $\lim_{n\rightarrow +\infty} \mathbb{E}(u_n)=+\infty$ by Lemma \ref{E decomposition lemma}.
\end{description}
\qed

\bigskip

(Xu Xu) School of Mathematics and Statistics, Wuhan University, Wuhan 430072, P.R. China

E-mail: xuxu2@whu.edu.cn\\[2pt]

(Chao Zheng) School of Mathematics and Statistics, Wuhan University, Wuhan 430072, P.R. China

E-mail: 2019202010023@whu.edu.cn\\[2pt]

\end{document}